\documentclass[12pt]{article}
\usepackage{latexsym,amsmath,amsthm}
\usepackage{amssymb}
\usepackage{graphicx}
\usepackage{appendix}
\usepackage{tikz}
\usepackage[colorlinks=true, linkcolor=blue]{hyperref}

\newtheorem{thm}{Theorem}[section]

\newtheorem{cor}[thm]{Corollary}
\newtheorem{lemma}[thm]{Lemma}
\newtheorem{conj}[thm]{Conjecture}
\newtheorem{claim}[thm]{Claim}

\theoremstyle{definition}
\newtheorem{prob}[thm]{Problem}

\theoremstyle{definition}

\theoremstyle{definition}
\newtheorem{defn}[thm]{Definition}

\theoremstyle{definition}

\theoremstyle{definition}

\theoremstyle{definition}

\theoremstyle{remark}
\newtheorem{conv}[thm]{Convention}

\theoremstyle{remark}

\usepackage[margin=1.2in]{geometry}

\usepackage{url}

\begin{document}

\title{Total non-negativity of some combinatorial matrices}

\author{David Galvin\thanks{Department of Mathematics,
University of Notre Dame, Notre Dame IN; dgalvin1@nd.edu. Supported in part by the Simons Foundation.}, Adrian Pacurar\thanks{Notre Dame, IN.}}

\maketitle

\begin{abstract}
Many combinatorial matrices --- such as those of binomial coefficients, Stirling numbers of both kinds, and Lah numbers --- are known to be totally non-negative, meaning that all minors (determinants of square submatrices) are non-negative.

The examples noted above can be placed in a common framework: for each one there is a non-decreasing sequence $(a_1, a_2, \ldots)$, and a sequence $(e_1, e_2, \ldots)$, such that the $(m,k)$ entry of the matrix is the coefficient of the polynomial $(x-a_1)\cdots(x-a_k)$ in the expansion of $(x-e_1)\cdots(x-e_m)$ as a linear combination of the polynomials $1, x-a_1, \ldots, (x-a_1)\cdots(x-a_m)$.   

We consider this general framework. For a non-decreasing sequence $(a_1, a_2, \ldots)$ we establish necessary and sufficient conditions on the sequence $(e_1, e_2, \ldots)$ for the corresponding matrix to be totally non-negative. As corollaries we obtain total non-negativity of matrices of rook numbers of Ferrers boards, and of graph Stirling numbers of chordal graphs.
\end{abstract}

\section{Introduction}

A matrix --- finite or infinite --- is {\em totally non-negative} if all minors (determinants of square sub-matrices) are non-negative. Totally non-negative matrices occur frequently in combinatorics and have been the subject of much investigation. See for example \cite{Brenti2, FominZelevinsky, GascaMicchelli, Skandera} for an overview.  
Here are a few of the most prominent examples:
\begin{itemize}
\item $\left[\binom{m}{k}\right]_{m,k \geq 0}$, where $\binom{m}{k}$ is the usual binomial coefficient;
\item $\left[{m \brace k}\right]_{m,k \geq 0}$, where ${m \brace k}$ is the {\em Stirling number of the second kind}, counting partitions of a set of size $m$ into $k$ non-empty blocks;
\item $\left[{m \brack k}\right]_{m,k \geq 0}$, where ${m \brack k}$ is the {\em (unsigned) Stirling number of the first kind}, counting partitions of a set of size $m$ into $k$ non-empty cyclically ordered blocks; and
\item $\left[L(m,k)\right]_{m,k \geq 0}$, where $L(m,k)$ is a {\em Lah number}, counting partitions of a set of size $m$ into $k$ non-empty linearly ordered blocks.
\end{itemize}
These examples can be placed in a common framework. Given two real sequences ${\bf a}=(a_1, a_2, \ldots)$ and ${\bf e}=(e_1, e_2, \ldots)$, either both infinite or both finite and of the same length, define a matrix $S^{{\bf a}, {\bf e}}=\left[S^{{\bf a}, {\bf e}}(m,k)\right]_{m,k \geq 0}$ via the relations
\begin{equation} \label{def-Sae}
\prod_{i=1}^m(x-e_i) = \sum_{k=0}^m S^{{\bf a}, {\bf e}}(m,k)\prod_{i=1}^k (x-a_i) 
\end{equation}
for $m \geq 0$. If ${\bf a}$ and ${\bf e}$ are infinite then $S^{{\bf a}, {\bf e}}$ is infinite with rows and columns indexed by $\{0, 1, 2, \ldots\}$, while if ${\bf a}$ and ${\bf e}$ are both of length $n$ then $S^{{\bf a}, {\bf e}}$ is $(n+1)$ by $(n+1)$ with rows and columns indexed by $\{0, 1, \ldots, n\}$. Note that (\ref{def-Sae}) uniquely determines $S^{{\bf a}, {\bf e}}(m,k)$ for each $m, k \geq 0$. In this framework,
\begin{itemize}
\item taking $e_i=-1$ and $a_i=0$ for all $i$ yields $S^{{\bf a}, {\bf e}}(m,k)=\binom{m}{k}$;
\item taking $e_i=0$ and $a_i=i-1$ for all $i$ yields $S^{{\bf a}, {\bf e}}(m,k)={m \brace k}$ via the identity 
\begin{equation} \label{eq-stir2-poly-def}
x^m = \sum_{k \geq 0}  {m \brace k}x(x-1)\cdots(x-(k-1))
\end{equation}
for $m \geq 0$; 
\item taking $e_i=-(i-1)$ and $a_i=0$ for all $i$ yields $S^{{\bf a}, {\bf e}}(m,k)={m \brack k}$ via the identity $x(x+1)\cdots(x+(m-1)) = \sum_{k \geq 0}  {m \brack k}x^k$ for $m \geq 0$; and
\item taking $e_i=-(i-1)$ and $a_i=i-1$ for all $i$ yields $S^{{\bf a}, {\bf e}}(m,k)=L(m,k)$ via the identity $x(x+1)\cdots(x+(m-1)) = \sum_{k \geq 0}  L(m,k)x(x-1)\cdots(x-(k-1))$
for $m \geq 0$. 
\end{itemize}

Another familiar object that fits into this framework is the collection of rook numbers of a Ferrers board. Let $b_1, b_2, \ldots$ be a non-decreasing sequence of non-negative integers, and let $B_m$ be the Ferrers board with $m$ columns that has $b_i$ cells in column $i$. The {\em rook number} $R_k(B_m)$ is the number of ways of placing $k$ non-attacking rooks on $B_m$. The factorization theorem of Goldman, Joichi and White \cite{GoldmanJoichiWhite2} says
$$
\sum_{k=0}^m R_{m-k}(B_m)x(x-1)\cdots (x-(k-1)) = \prod_{i=1}^m (x+b_i-i+1).
$$   
Taking $e_i = i-1-b_i$ and $a_i=i-1$, we see from (\ref{def-Sae}) that $S^{{\bf a}, {\bf e}}(m,k) = R_{m-k}(B_m)$.

The main result of this note is a characterization, for each non-decreasing sequence ${\bf a}$, of those sequences ${\bf e}$ such that the matrix $S^{{\bf a}, {\bf e}}$ is totally non-negative.
\begin{defn}
If ${\bf a}$ is non-decreasing, we say that ${\bf e}$ is a {\em restricted growth sequence} relative to ${\bf a}$ if for each $i \geq 1$ it holds that $e_i \leq a_{f(i)}$, where $f(1)=1$ and for $i \geq 1$
$$
f(i+1) = \left\{
\begin{array}{cl}
f(i) & \mbox{if $e_i < a_{f(i)}$} \\
f(i)+1 & \mbox{if $e_i = a_{f(i)}$}. 
\end{array}
\right.
$$  
\end{defn}
In other words, each $e_i$ is at most a certain cap. The cap for $e_1$ is $a_1$. If $e_1 < a_1$ then the cap for $e_2$ is also $a_1$, while if $e_1=a_1$ then the cap for $e_2$ is $a_2$. In general, the cap for $e_i$ is some $a_{i'}$, and if $e_i<a_{i'}$ then the cap for $e_{i+1}$ is also $a_{i'}$, while if $e_i=a_{i'}$ then the cap for $e_{i+1}$ is $a_{i'+1}$. If ${\bf a}=(0,1,\ldots, n-1, \ldots)$ then a non-negative integral sequence ${\bf e}$ is a restricted growth sequence relative to ${\bf a}$ exactly if it is a restricted growth sequence in the usual sense, that is, one satisfying $e_1=0$ and $e_{i+1} \leq 1+\max_{j=1, \ldots, i} e_j$ for $i \geq 1$.

Notice that in the examples of binomial coefficients, Stirling numbers of both kinds and Lah numbers above, ${\bf a}$ is non-decreasing and ${\bf e}$ is a restricted growth sequence relative to ${\bf a}$. The total non-negativity of the matrices arising from these examples is thus recovered from the following result.
\begin{thm} \label{thm-main}
Let ${\bf a}$ be a non-decreasing sequence. Then
\begin{enumerate}
\item \label{mainthm-item1} the matrix $S^{{\bf a}, {\bf e}}$ is totally non-negative if and only if ${\bf e}$ is a restricted growth sequence relative to ${\bf a}$, and
\item \label{mainthm-item2} if ${\bf e}$ is {\em not} a restricted growth sequence relative to ${\bf a}$ then the failure of $S^{{\bf a}, {\bf e}}$ to be totally non-negative is witnessed by a negative entry in $S^{{\bf a}, {\bf e}}$.
\end{enumerate}
\end{thm}
Also, since the sequence ${\bf a}=(0,1,2,\ldots)$ is non-decreasing and for non-negative $(b_1, b_2, \ldots)$  the sequence ${\bf e}=(-b_1, 1-b_2, 2-b_3)$ is restricted growth relative to ${\bf a}$, we immediately have the following corollary concerning matrices of rook numbers.
\begin{cor} \label{cor-rook}
If $(b_1, b_2, \ldots)$ is a non-deceasing sequence of positive integers, and $B_m$ is the Ferrers board with $m$ columns and with $b_i$ cells in the $i$th column, then the matrix $\left(R_{m-k}(B_m)\right)_{m, k \geq 0}$ is totally non-negative, where $R_k(B_m)$ is the number of ways of placing $k$ non-attacking rooks on $B_m$. 
\end{cor} 

En route to proving Theorem \ref{thm-main} we will show that for arbitrary ${\bf a}$ and ${\bf e}$ total non-negativity of $S^{{\bf a}, {\bf e}}$ is guaranteed by the condition $\inf {\bf a} \geq \sup {\bf e}$ (equivalently, $a_i-e_j \geq 0$ for all $i, j$). This observation (Corollary \ref{cor-obsv-Gon}, suggested to us by Gonzales \cite{Gonzales}) represents all we can say at present when ${\bf a}$ is not assumed to be non-decreasing. 
\begin{prob}
For ${\bf a}$ which is not non-decreasing, characterize those ${\bf e}$ for which $S^{{\bf a}, {\bf e}}$ is totally non-negative.
\end{prob} 

The proof of Theorem \ref{thm-main} involves producing a weighted planar network whose path matrix is $S^{{\bf a}, {\bf e}}$, and then appealing to Lindstr\"om's lemma (see Section \ref{sec-proof-of-main-theorem} for details).  The network that we initially produce, however, does not have all non-negative entries, precluding an immediate application of Lindstr\"om. A substantial part of the proof involves carefully modifying the weights of the initial network to remove the negative entries, without changing the associated path matrix. We prove Theorem \ref{thm-main} in Section \ref{sec-proof-of-main-theorem}. Before that, in Section \ref{sec-chordal-app}, we consider an application to graph Stirling numbers of chordal graphs.

\medskip

The numbers $S^{{\bf a}, {\bf e}}(m,k)$ defined in (\ref{def-Sae}) satisfy the recurrence
\begin{equation} \label{Sae-rec}
S^{{\bf a}, {\bf e}}(m,k) = S^{{\bf a}, {\bf e}}(m-1,k-1) + (a_{k+1}-e_m)S^{{\bf a}, {\bf e}}(m-1,k) ~~~\mbox{for $m,k > 0$}
\end{equation}
with initial conditions $S^{{\bf a}, {\bf e}}(0,0) = 1$, $S^{{\bf a}, {\bf e}}(0,k)  = 0$ for $k > 0$ and $S^{{\bf a}, {\bf e}}(m,0) = \prod_{i=1}^m (a_1-e_i)$ for $m > 0$ (we prove this in Section \ref{sec-proof-of-main-theorem}, see (\ref{recurrence})). Various forms of this recurrence have appeared in the literature. As observed in \cite{Gonzales},
with suitable choices of ${\bf a}$ and ${\bf e}$ the recurrence (\ref{Sae-rec}) can encode
\begin{itemize}
\item some generalizations of the classical rook numbers \cite{CelesteCorcinoGonzales},
\item the normal order coefficients of the word $(VU)^n$ in the Weyl algebra generated by symbols $V, U$ satisfying $UV - V U = hV^s$ \cite{CelesteCorcinoGonzales},
\item Hsu and Shiue's generalized Stirling numbers \cite{HsuShiue}, 
\item the Jacobi-Stirling numbers (coefficients of the Jacobi differential operator) \cite{CelesteCorcinoGonzales, EverittKwonLittlejohnWellmanYoon},
\end{itemize}
as well as encoding Binomial coefficients, Stirling numbers of both kinds, Lah numbers and rook numbers.

A number of authors have considered the question of total non-negativity of matrices $[a_{m,k}]_{m, k \geq 0}$ with the $a_{m,k}$ defined via recurrences similar to (\ref{Sae-rec}). Brenti \cite{Brenti2}, for example, considered the recurrence $a_{m,k}= z_ma_{m-t,k-1}+y_ma_{m-1,k-1}+x_ma_{m-1,k}$ ($t \in {\mathbb N}$). More recently Chen, Liang and Wang \cite{ChenLiangWang2, ChenLiangWang}
 considered $a_{m,k}=r_ka_{m-1,k-1} + s_ka_{m-1,k} + t_{k+1}a_{m-1,k+1}$ and also the more general situation where the $a_{m,k}$'s form a Riordan array.
The recurrence (\ref{Sae-rec}) does not seem to fit these settings. 

To conclude the introduction, we mention a nice conjecture of Brenti \cite[Conjecture 6.10]{Brenti3} to which the present work may be related. The {\em Eulerian number} $A(m, k)$ is the number of permutations of $\{1,\ldots, m\}$ with exactly $k$ ascents. It satisfies a recurrence that is very similar to (\ref{Sae-rec}), namely
$$
A(m,k) = (m-k)A(m-1,k-1)+(k+1)A(m-1,k).
$$
\begin{conj}
The matrix $\left[A(m,k)\right]_{m,k \geq 0}$ is totally non-negative.
\end{conj}

\section{Graph Stirling numbers of chordal graphs} \label{sec-chordal-app}

The Stirling numbers of the second kind have a natural generalization to the setting of graphs. For a graph $G$ and an integer $k$, the {\em graph Stirling number of the second kind} ${G \brace k}$ is the number of ways of partitioning the vertex set of $G$ into $k$ non-empty independent sets (an {\em independent set} being a set of pairwise non-adjacent vertices). This is indeed a generalization, since if $E_m$ is the graph on $m$ vertices with no edges, then ${E_m \brace k}={m \brace k}$. 

This notion of graph Stirling number of the second kind was probably first introduced by Tomescu \cite{Tomescu} and was subsequently reintroduced by numerous authors including Korfhage \cite{Korfhage}, Goldman, Joichi and White \cite{GoldmanJoichiWhite} and Duncan and Peele \cite{DuncanPeele}. Its properties have been well studied, see for example \cite{Brenti, BrentiRoyleWagner, DoGalvin, EngbersGalvinHilyard, FarrellWhitehead, MohrPorter, Munagi}.

The Stirling number of the first kind does not have such a natural graph analog. In \cite{EuFuLiangWong} Eu, Fu, Liang and Wong present a notion of a graph Stirling number of the first kind for the family of quasi-threshold graphs, based on generalizations of the relation $x^mD^m = \sum_{k \geq 0} (-1)^{m-k}{m \brack k} (xD)^k$ in the Weyl algebra on symbols $x$ and $D$ (the algebra over the reals generated by the relation $Dx = xD + 1$).  

Here we take a different approach. It is well-known that the inverse of the matrix of Stirling numbers of the second kind is the matrix of {\em signed} Stirling numbers of the first kind:
$$
\left[{m \brace k}\right]_{m,k\geq 0}^{-1} = \left[(-1)^{m-k}{m \brack k}\right]_{m,k\geq 0}.
$$
This suggests the following. For a graph $G$ on $n$ vertices, ordered $v_1, \ldots, v_n$, let $G_m$ denote the subgraph of $G$ induced by $v_1, \ldots, v_m$, and consider the matrices
$$
S_G = \left[{G_m \brace k}\right]_{m,k =0}^n~~~\mbox{and}~~~s_G = S_G^{-1}. 
$$  
So the $(m,k)$ entry of $s_{E_n}$ is $(-1)^{m-k}{m \brack k}$, and the (non-negative) quantity ${m \brack k}$ has a clean combinatorial interpretation, as the size of a set of permutations. 

It would be of interest to have a combinatorial interpretation of the absolute value of the $(m,k)$ entry of $s_G$ for general $G$, leading to a combinatorial notion of graph Stirling numbers of the first kind for all graphs. To this end, it would be helpful to know the sign of the $(m,k)$ entry of $s_G$. The reason for this is as follows. If $M$ is a lower triangular matrix with non-negative integer entries and $1$'s down the diagonal (note $S_G$ is of this form) then it is possible to express the $(m,k)$ entry of $M^{-1}$ as $\sum_{a \in {\mathcal A}(m,k)} {\rm sign}_{m,k}(a)$ where ${\mathcal A}(m,k)$ is some combinatorially defined set and ${\rm sign}_{m,k}$ is a sign function taking values in $\{1,-1\}$ (see, for example, \cite{BapatGhorbani}, where ${\mathcal A}(m,k)$ is a certain set of paths in a bipartite multigraph).  If the sign of the $(m,k)$ entry of $M^{-1}$ is known, and happens to be positive, then one could find a combinatorial interpretation of the entry as a count of a set (rather than a signed count of a set, or as the difference in the sizes of two sets) by constructing an injection from ${\mathcal A}^-:=\{a \in {\mathcal A}(m,k):{\rm sign}_{m,k}=-1\}$ into ${\mathcal A}^+:=\{a \in {\mathcal A}(m,k):{\rm sign}_{m,k}=1\}$, and finding a description of those $a \in {\mathcal A}^+$ that are not in the range of the injection. If the sign of the $(m,k)$ entry of $M^{-1}$ is known to be negative, one would seek instead an injection from ${\mathcal A}^+$ into ${\mathcal A}^-$. See, for example, \cite{EngbersGalvinSmyth}, where this strategy is employed to provide combinatorial interpretations of entries of inverses of matrices of certain restricted Stirling and Lah numbers. 

It is easy to find examples of graphs $G$ such that however the vertices are ordered the pattern of signs in the matrix $s_G$ is quite chaotic, making the approach just discussed difficult to implement. However, there is a class of graphs which admit a natural ordering of the vertices with respect to which the pattern of signs in $s_G$ is very well behaved, and in fact  has the same checkerboard sign pattern as $s_{E_n}$, that is, with the $(m,k)$ entry having sign $(-1)^{m-k}$.  A {\em chordal graph} is a graph in which every cycle of length four or greater has a chord, that is, it is a graph that contains no induced cycles of length four or greater. A useful characterization of chordal graphs is that $G$ is chordal if and only if it is possible to order the vertices as $v_1, \ldots, v_n$ so that for each $m \in \{1,\ldots, n\}$ the neighbors of $v_m$ among $v_1, \ldots, v_{m-1}$ induce a clique (see for example \cite[Section 5.3]{West}). Such an ordering is referred to as a {\em perfect elimination order}.
\begin{thm} \label{thm-chordal}
Let $G$ be a chordal graph with perfect elimination order $v_1, \ldots, v_n$, and let $G_m$ be the subgraph of $G$ induced by $v_1, \ldots, v_m$. Let $S_G = \left[{G_m \brace k}\right]_{m,k =0}^n$ and $s_G = S_G^{-1}$. For all $m, k$  the $(m,k)$ entry of $s_G$ has sign $(-1)^{m-k}$.
\end{thm}   

A stronger result than Theorem \ref{thm-chordal} holds. Notice that the matrix $S_G$ has determinant $1$ and so by Cramer's rule the $(m,k)$ entry of the inverse is $(-1)^{m-k}$ times the determinant of the $n-1$ by $n-1$ minor obtained from $S_G$ by deleting the $k$th row and the $m$th column. It follows that if $S_G$ is totally non-negative then the $(m,k)$ entry of $s_G$ has sign $(-1)^{m-k}$, and so the following result generalizes Theorem \ref{thm-chordal}.
\begin{thm} \label{thm-totallynn}
Let $G$ be a chordal graph with perfect elimination order $v_1, \ldots, v_n$, and let $G_m$ be the subgraph of $G$ induced by $v_1, \ldots, v_m$. Let $S_G = \left[{G_m \brace k}\right]_{m,k =0}^n$. Then $S_G$ is totally non-negative.
\end{thm}    

As we will now see, Theorem \ref{thm-totallynn} is a special case of Theorem \ref{thm-main}. The {\em chromatic polynomial} $\chi_G(x)$ of a graph $G$ is the polynomial in $x$ whose value at positive integers $x$ is the number of ways of coloring $G$ from a palette of $x$ colors in such a way that adjacent vertices receive distinct colors. That $\chi_G(x)$ is indeed a polynomial in $x$ follows from the following identity: for $G$ a graph on $m$ vertices,
\begin{equation} \label{eq-chromatic_polynomial}
\chi_G(x) = \sum_{k=0}^m {G \brace k} x(x-1)\cdots(x-(k-1)).
\end{equation}
Indeed, one way to enumerate the colorings of $G$ from a palette of $x$ colors in such a way that adjacent vertices receive distinct colors is to first specify $k$, the number of colors used, then specify a partition of the vertex set of $G$ into $k$ non-empty independent sets (${G \brace k}$ options), which will be the color classes, and finally specify the colors that appear on each of the classes ($x(x-1)\cdots(x-(k-1))$ options). (Observe that by taking $G$ to be the graph on $m$ vertices with no edges we recover (\ref{eq-stir2-poly-def}) from (\ref{eq-chromatic_polynomial})).

For a chordal graph $G$ with perfect elimination order $v_1, \ldots, v_n$, for $i\geq 1$ denote by $e_i=e_i(G)$ the number of neighbors that $v_i$ has among $v_1, \ldots, v_{i-1}$. We have that   $\chi_{G_m}(x) = (x-e_1)(x-e_2)\cdots(x-e_m)$ (coloring the vertices of $G_m$ sequentially from $v_1$ to $v_m$, at the step when $v_j$ is colored all colors are available except those used on the $e_j$ neighbors of $v_j$ among $\{v_1, \ldots, v_{j-1}\}$; since these neighbors form a clique, between them they account for $e_j$ colors, leaving $x-e_j$ available for $v_j$). Thus, in light of (\ref{eq-chromatic_polynomial}), if we knew that $(e_1, \ldots, e_n)$ formed a restricted growth sequence relative to $(0,1,\ldots, n-1)$, then the total non-negativity of $\left[{G_m \brace k}\right]_{m,k=0}^n$ would follow from Theorem \ref{thm-main}.

In fact, we have the following.
\begin{claim} \label{clm-restricted-growth-chordal}
Let $G$ be a chordal graph with perfect elimination order $v_1, \ldots, v_n$. Defining $e_i=e_i(G)$ as above, we have that $(e_1, \ldots, e_n)$ is a restricted growth sequence relative to $(0,1,\ldots, n-1)$. Moreover, if $(e'_1, \ldots, e'_n)$ is any restricted growth sequence relative to $(0,1,\ldots, n-1)$ then there is a chordal graph $G$ with perfect elimination order $v_1, \ldots, v_n$ such that $e_i(G)=e'_i$ for all $i\leq n$.
\end{claim}

\medskip

\noindent {\em Proof}:
We begin by showing that $(e_1, \ldots, e_n)$ is a restricted growth sequence relative to $(0,1,\ldots, n-1)$. Certainly $e_1=0$. Now consider vertex $v_k$ for $k > 1$. It is adjacent to $e_k$ vertices among $v_1, \ldots, v_{k-1}$, with the largest of these (in the ordering $v_1 < v_2 < \cdots$) being, say, $v_j$. Because $v_k$ forms a clique with its neighbors among $v_1, \ldots, v_{k-1}$, it follows that $v_j$ has at least $e_k-1$ neighbors among $v_1, \ldots, v_{j-1}$, so $e_j \geq e_k-1$. From this it follows that $e_k \leq e_j+1 \leq 1+\max_{i < k} e_i$, exactly the condition that says that $(e_1, \ldots, e_n)$ is a restricted growth sequence relative to $(0,1,\ldots, n-1)$.     

Next suppose $(e'_1, \ldots, e'_n)$ is a restricted growth sequence relative to $(0,1,\ldots, n-1)$. We inductively construct a chordal graph $G$ with perfect elimination order $v_1, \ldots, v_n$ such that $e_i(G)=e'_i$ for all $i=1,\ldots, n$, starting with an isolated vertex $v_1$. Suppose that the adjacency structure among $v_1, \ldots, v_{k-1}$ has been determined. We have that $e_k \leq 1+\max_{i < k} e_i$, which means that (by induction) among $v_1, \ldots, v_{k-1}$ there are some $e_k$ vertices that form a clique. The construction can be continued by joining $v_k$ to any such $e_k$ vertices.
\qed

\section{Proof of Theorem \ref{thm-main}} \label{sec-proof-of-main-theorem}

A key tool will be the following explicit expression for $S^{{\bf a}, {\bf e}}(m,k)$.
\begin{lemma} \label{lem-gabes-formula}
For arbitrary ${\bf a}$ and ${\bf e}$ we have
\begin{equation} \label{formula-Gabe}
S^{{\bf a}, {\bf e}}(m,k) = \sum_{\scriptstyle S=\{s_1,\ldots,s_{m-k}\} \subseteq \{1,\ldots,m\} \atop \scriptstyle s_1<\ldots<s_{m-k}} \prod_{i=1}^{m-k} (a_{s_i-i+1} -e_{s_i}),
\end{equation}
and, equivalently, denoting the $(m,k)$ entry of $(S^{{\bf a}, {\bf e}})^{-1}$ by $s^{{\bf a}, {\bf e}}(m,k)$,
\begin{equation} \label{formula-Gabe2}
(-1)^{m-k}s^{{\bf a}, {\bf e}}(m,k) = \sum_{\scriptstyle S=\{s_1,\ldots,s_{m-k}\} \subseteq \{1,\ldots,m\} \atop \scriptstyle s_1<\ldots<s_{m-k}} \prod_{i=1}^{m-k} (a_{s_i} - e_{s_i-i+1}).
\end{equation}
\end{lemma}

Notice that in the chordal graph setting $e_{s_i-i+1}$ is the number of edges from vertex $v_{s_i-i+1}$ to earlier vertices, so is at most $s_i-i$, which is at most $s_i-1$, which is $a_{s_i}$, and so the quantity on the right-hand side of the formula for $(-1)^{m-k}s(m,k)$ is non-negative. This establishes directly that the sign of the $(m,k)$ entry of $s_G$ is $(-1)^{m-k}$, as asserted by Theorem \ref{thm-chordal}.

\medskip

\noindent {\em Proof (of Lemma \ref{lem-gabes-formula})}: 
We begin by noting that (\ref{formula-Gabe}) implies (\ref{formula-Gabe2}). Indeed, from (\ref{def-Sae}) we have that $S^{{\bf a}, {\bf e}}(m,k)$ is the coefficient of $(x-a_1)\cdots(x-a_k)$ in the unique expansion of $(x-e_1)\cdots(x-e_m)$ as a linear combination of $1, x-a_1, \ldots,  (x-a_1)\cdots(x-a_m)$, and so from basic linear algebra considerations we see that the $s^{{\bf a}, {\bf e}}(m,k)$ are uniquely determined by the relations
$$
\prod_{i=1}^m(x-a_i) = \sum_{k=0}^m s^{{\bf a}, {\bf e}}(m,k)\prod_{i=1}^k(x-e_i)
$$
for $m \geq 0$. Since in (\ref{formula-Gabe}) ${\bf a}$ and ${\bf e}$ are arbitrary, a direct application of that identity yields
\begin{eqnarray*}
(-1)^{m-k}s^{{\bf a}, {\bf e}}(m,k) & = & (-1)^{m-k}\sum_{\scriptstyle S=\{s_1,\ldots,s_{m-k}\} \subseteq \{1,\ldots,m\} \atop \scriptstyle s_1<\ldots<s_{m-k}} \prod_{i=1}^{m-k} (e_{s_i-i+1} -a_{s_i}) \\
& = & \sum_{\scriptstyle S=\{s_1,\ldots,s_{m-k}\} \subseteq \{1,\ldots,m\} \atop \scriptstyle s_1<\ldots<s_{m-k}} \prod_{i=1}^{m-k} (a_{s_i} - e_{s_i-i+1}).
\end{eqnarray*}
Of course the same argument in reverse shows that also (\ref{formula-Gabe2}) implies (\ref{formula-Gabe}).

We now show that both sides of (\ref{formula-Gabe}) satisfy the same recurrence relation and initial conditions. To that end write $f(m,k)$ for the expression on the right-hand side of (\ref{formula-Gabe}). We begin by establishing some boundary values for $f(m,k)$.
\begin{itemize}
\item We have $f(0,0)=1$ (the sum has one summand, associated with $S=\emptyset$, and this summand is the empty product and so has value $1$), and more generally $f(m,m)=1$ for all $m$.
\item For $m>0$, $f(m,0)=(a_1-e_1)\cdots (a_1-e_m)$.
\item For $k>0$, $f(0,k)=0$ (the sum defining $f$ in this case is empty), and more generally for $k>m$, $f(m,k)=0$.
\end{itemize}
Next we establish a recurrence for $f(m,k)$. For $m > k > 0$ we have 
$$
f(m,k) = f(m-1,k-1) + (a_{k+1}-e_m)f(m-1,k).
$$
The terms on the right-hand side here come from considering first those $S$ with $m \not \in S$ and then those with $m \in S$; in the latter case $m$ is always the greatest element of $S$ and so contributes a factor $a_{m-(m-k)+1}-e_m=a_{k+1}-e_m$ to each summand.

Next consider the quantity $S^{{\bf a}, {\bf e}}(m,k)$. We easily have $S^{{\bf a}, {\bf e}}(0,0)=1$, and more generally $S^{{\bf a}, {\bf e}}(m,m)=1$ for all $m$, as well as 
$S^{{\bf a}, {\bf e}}(m,0)= (a_1-e_1)\cdots (a_1-e_m)$ for $m>0$ (evaluate both sides of (\ref{def-Sae}) at $x=a_1$). We also have $S^{{\bf a}, {\bf e}}(0,k)=0$ for $k > 0$ and more generally $S^{{\bf a}, {\bf e}}(m,k)=0$ for $k > m$. We also have the recurrence   
\begin{equation} \label{recurrence}
S^{{\bf a}, {\bf e}}(m,k) = S^{{\bf a}, {\bf e}}(m-1,k-1) + (a_{k+1}-e_m) S^{{\bf a}, {\bf e}}(m-1,k)
\end{equation}
for $m > k > 0$. To verify this, consider the expression
\begin{equation} \label{local-1}
\begin{array}{c}
S^{{\bf a}, {\bf e}}(m,0) + S^{{\bf a}, {\bf e}}(m-1,m-1)(x-a_1)\cdots (x-a_m) \\
+\sum_{k=1}^{m-1} \left(S^{{\bf a}, {\bf e}}(m-1,k-1)+(a_{k+1}-e_m)S^{{\bf a}, {\bf e}}(m-1,k)\right)(x-a_1)\cdots(x-a_k)
\end{array}
\end{equation}
(a linear combination of the polynomials $1, x-a_1, \ldots, (x-a_1)\cdots (x-a_m)$).
Rearranging terms (\ref{local-1}) becomes
\begin{equation} \label{local-2}
\begin{array}{c}
S^{{\bf a}, {\bf e}}(m,0) +  S^{{\bf a}, {\bf e}}(m-1,0)(x-a_1)\\
+(x-e_m)\sum_{k=1}^{m-1} S^{{\bf a}, {\bf e}}(m-1,k)(x-a_1)\cdots(x-a_k).
\end{array}
\end{equation}
Writing $x-a_1 = (x-e_m)-(a_1-e_m)$ in the second term of (\ref{local-2}) yields
\begin{equation} \label{local-3}
\begin{array}{c}
S^{{\bf a}, {\bf e}}(m,0) - S^{{\bf a}, {\bf e}}(m-1,0)(a_1-e_m)\\
+(x-e_m)\sum_{k=0}^{m-1} S^{{\bf a}, {\bf e}}(m-1,k)(x-a_1)\cdots(x-a_k).
\end{array}
\end{equation}
Via the initial conditions the first two terms of (\ref{local-3}) sum to $0$, and via the defining relation for $S^{{\bf a}, {\bf e}}(m-1,\cdot)$ ((\ref{def-Sae}) with $m$ replaced by $m-1$) the remaining terms sum to $\prod_{i=1}^m (x-e_i)$. The recurrence (\ref{recurrence}) now follows from (\ref{def-Sae}) via linear algebra considerations.

Since $f(m,k)$ and $S^{{\bf a}, {\bf e}}(m,k)$ satisfy the same recurrence and initial conditions, they are equal.
\qed

\medskip

Lemma \ref{lem-gabes-formula} allows us to express $S^{{\bf a}, {\bf e}}(m,k)$ in terms of complete symmetric and elementary symmetric functions. Denote by $h_{\ell}(x_1,\ldots, x_t)$ the degree $\ell$ complete symmetric polynomial in $x_1,\ldots, x_t$ (the sum of all degree $\ell$ monomials with coefficients $1$) and by $s_{\ell}(x_1,\ldots, x_t)$ the degree $\ell$ elementary symmetric polynomial in $x_1,\ldots, x_t$ (the sum of all degree $\ell$ linear monomials with coefficients $1$); so, for example, 
$h_2(x_1,x_2,x_3) = x_1^2+x_2^2+x_3^2 + x_1x_2 + x_1x_3+x_2x_3$
while
$s_2(x_1,x_2,x_3) = x_1x_2 + x_1x_3+x_2x_3$. (We will later also use $s_{\cdot}$ for a source in a planar network, but this should not cause confusion as the meaning will be clear from the context.)

\begin{lemma} \label{lem-sym}   
For arbitrary ${\bf a}$ and ${\bf e}$, 
\begin{equation} \label{sym-polys}
S^{{\bf a}, {\bf e}}(m,k) = \sum_{\ell=0}^{m-k} (-1)^\ell h_{m-k-\ell}(a_1,\ldots,a_{k+1})s_\ell(e_1,\ldots,e_m). 
\end{equation}
\end{lemma}

\medskip

\noindent {\em Proof}:  One possible approach is to show that that the expressions on the right-hand sides of (\ref{formula-Gabe}) and (\ref{sym-polys}) are equal. This can be achieved by noting that when the right-hand side of (\ref{formula-Gabe}) is expanded as a polynomial in the $e_i$'s, the monomials that arise are precisely the linear monomials in $e_1, \ldots, e_m$. For a given $\ell$, $0 \leq \ell \leq m-k$, and $T=\{t_1, \ldots, t_\ell\} \subseteq \{1,\ldots,m\}$ with $t_1<\ldots < t_\ell$, the coefficient of $e_{t_1}\ldots e_{t_\ell}$ turns out to be $(-1)^\ell h_{m-k-\ell}(a_1,\ldots,a_{k+1})$ (independent of the particular choice of $T$); this proves the lemma.

We take instead a linear algebra approach. From (\ref{def-Sae}) we have
$$
\prod_{i=1}^m (x-e_i) = \sum_k S^{{\bf a}, {\bf e}}(m,k)\prod_{i=1}^k (x-a_i)
$$
(where the sum runs over all integers $k$, although the summand will only be non-zero for $k \in \{0,1, \ldots, m\}$). It follows that
\begin{eqnarray*}
\prod_{i=1}^m(x-e_i) & = & \sum_j S^{{\bf 0}, {\bf e}}(m,j)x^j \\
& = & \sum_j S^{{\bf 0}, {\bf e}}(m,j) \sum_k S^{{\bf a}, {\bf 0}}(j,k) \prod_{i=1}^k (x-a_i) \\
& = & \sum_k \left(\sum_j S^{{\bf 0}, {\bf e}}(m,j) S^{{\bf a}, {\bf 0}}(j,k)\right) \prod_{i=1}^k (x-a_i)
\end{eqnarray*}
so that
\begin{equation} \label{newint1}
S^{{\bf a}, {\bf e}}(m,k) = \sum_j S^{{\bf 0}, {\bf e}}(m,j) S^{{\bf a}, {\bf 0}}(j,k).
\end{equation}
Now  Lemma \ref{lem-gabes-formula} gives
\begin{eqnarray}
S^{{\bf 0}, {\bf e}}(m,j) & =  & \sum_{\scriptstyle S=\{s_1,\ldots,s_{m-j}\} \subseteq \{1,\ldots,m\} \atop \scriptstyle s_1<\ldots<s_{m-j}} \prod_{i=1}^{m-j} (-e_{s_i}) \nonumber \\
& = & (-1)^{m-j}s_{m-j}(e_1,\ldots, e_m) \label{comp1}
\end{eqnarray}
and
\begin{eqnarray}
S^{{\bf a}, {\bf 0}}(j,k) & = & \sum_{\scriptstyle S=\{s_1,\ldots,s_{j-k}\} \subseteq \{1,\ldots,j\} \atop \scriptstyle s_1<\ldots<s_{j-k}} \prod_{i=1}^{j-k} a_{s_i-i+1}  \nonumber \\
& = & h_{j-k}(a_1,\ldots, a_{k+1}). \label{comp2}
\end{eqnarray}
Combining (\ref{comp1}) and (\ref{comp2}) with (\ref{newint1}), and re-indexing via $\ell=m-j$, leads to (\ref{sym-polys}).
\qed

\medskip

We now require some well-known results from the theory of totally non-negative matrices. 
A {\em planar network} $P$ is a directed planar graph with a subset of vertices designated as sources and a subset of vertices designated as sinks. A {\em weighted} planar network $(P,w)$ is a planar network $P$ together with a function $w:\vec{E}(P)\rightarrow {\mathbb R}$, which we think of as an assignment of weights to the edges of $P$. Figure \ref{fig1} shows a particular weighted planar network with sources $\{s_i:i=0,1,2,\ldots\}$ and sinks $\{t_i:i=0,1,2,\ldots\}$. The other vertices are the points of intersection between the vertical and horizontal lines. Horizontal lines are oriented to the right and vertical lines are oriented upward. All horizontal edge weights are $1$, while the weights of the vertical edges are given by the $x_{ij}$'s. 

\begin{figure}[ht!]
\begin{center}
\begin{tikzpicture}
\draw [thick, ->] (-6,6) -- (6,6);
\draw [thick, ->] (-6,5) -- (6,5);
\draw [thick, ->] (-6,4) -- (6,4);
\draw [thick, ->] (-6,3) -- (6,3);
\draw [thick, ->] (-6,1) -- (6,1);
\draw [thick, ->] (-6,0) -- (6,0);
\draw [thick, ->] (-5,-.8) -- (-5,6);
\draw [thick, ->] (-3.5,-.8) -- (-3.5,5);
\draw [thick, ->] (-2,-.8) -- (-2,4);
\draw [thick, ->] (-.5,-.8) -- (-.5,3);
\draw [thick, ->] (5,-.8) -- (5,1); 
\node [left] at (-6,6) {$s_0$};
\node [left] at (-6,5) {$s_1$};
\node [left] at (-6,4) {$s_2$};
\node [left] at (-6,3) {$s_3$};
\node [left] at (-6,2) {$\vdots$};
\node [left] at (-6,1) {$s_{n-1}$};
\node [left] at (-6,0) {$s_n$};
\node [right] at (6,6) {$t_0$};
\node [right] at (6,5) {$t_1$};
\node [right] at (6,4) {$t_2$};
\node [right] at (6,3) {$t_3$};
\node [right] at (6,2) {$\vdots$};
\node [right] at (6,1) {$t_{n-1}$};
\node [right] at (6,0) {$t_n$};
\node [right] at (-5,5.5) {$\!x_{11}$};
\node [right] at (-5,4.5) {$\!x_{21}$};
\node [right] at (-5,3.5) {$\!x_{31}$};
\node [right] at (-5,.5) {$\!x_{n1}$};
\node [right] at (-3.5,4.5) {$\!x_{22}$};
\node [right] at (-3.5,3.5) {$\!x_{32}$};
\node [right] at (-3.5,.5) {$\!x_{n2}$};
\node [right] at (-2,3.5) {$\!x_{33}$};
\node [right] at (-2,.5) {$\!x_{n3}$};
\node [right] at (-.5,.5) {$\!x_{n4}$};
\node [right] at (5,.5) {$\!x_{nn}$};
\node at (1,2) {$\ddots$};
\node [left] at (-6,-1) {$\vdots$};
\node at (1,-1) {$\vdots$};
\node [right] at (6,-1) {$\vdots$};
\end{tikzpicture}
\caption{A weighted planar network.} \label{fig1}
\end{center}
\end{figure}

The planar network shown in Figure \ref{fig1} is the only one that we will consider in the sequel, and while we will consider many weight functions, they will all have the same form as that shown in Figure \ref{fig1} (that is, with all horizontal edge weights being $1$). We will represent a generic such weight function by a doubly infinite lower triangular array, viz:
\begin{equation} \label{array-of-weights}
\begin{array}{cccccccc}
x_{11} & & & & & & \\
x_{21} & x_{22} & & & & & \\
x_{31} & x_{32} & x_{33} & & & & \\
x_{41} & x_{42} & x_{43} & x_{44} & & & \\
x_{51} & x_{52} & x_{53} & x_{54} & x_{55} &   & \\
\vdots & \vdots & \vdots & \vdots &        & \ddots & \\
x_{m1} & x_{m2} & x_{m3} & x_{m4} & \cdots & x_{m(m-1)} & x_{mm} \\
\vdots & \vdots & \vdots & \vdots &  & & & \ddots.
\end{array}
\end{equation}
We will refer to the location of the weight $x_{mk}$ in this array as the {\em $[m,k]$ position} of the array (using square brackets to distinguish this from an entry in a matrix), and we will refer to the corresponding edge of the planar network as the {\em $[m,k]$ edge}. For each fixed $k$, the collection of $[m,k]$ edges, $m=k, k+1, \ldots$, is referred to as the {\em $k$-th column} of the planar network.

By the {\em path matrix} of a weighted planar network we mean the doubly infinite matrix whose $(i,j)$ entry (with rows and columns indexed by $\{0,1,2,\ldots\}$) is the sum of the weights of all the directed paths from $s_i$ to $t_j$, where the weight of one such path is the product, over all edges traversed, of the weight of the edge. For example, the $(3,1)$ entry of the path matrix of the weighted planar network shown in Figure \ref{fig1} is $x_{31}x_{21}+ x_{31}x_{22}+ x_{32}x_{22}$. Notice that the path matrix of the weighted planar network shown in Figure \ref{fig1} is lower-triangular with $1$'s down the main diagonal. The following \cite{Lindstrom} (see also, for example, \cite{Skandera}) is a standard result from the theory of totally non-negative matrices.
\begin{lemma} (Lindstr\"om's Lemma) \label{lem-tp}
If the matrix $M$ is the path matrix of a weighted planar network in which all weights are non-negative, then $M$ is totally non-negative. 
\end{lemma}
Lindstr\"om's Lemma in fact says more: the minor corresponding to selecting the rows indexed by $I$ and columns indexed by $J$ (with indexing of rows and columns starting from $0$) equals the sum of the weights of all the collections of $|I|$ vertex disjoint paths from the sources $\{s_i: i \in I\}$ to the sinks $\{t_j: j \in J\}$, where the weight of a collection of paths is the product of the weights of the individual paths in the collection. We will not need this level of precision in our analysis.

In the presence of Lindstr\"om's Lemma, to prove item \ref{mainthm-item1} of Theorem \ref{thm-main} it suffices to construct, for each ${\bf a}$ and ${\bf e}$ with ${\bf a}$ non-decreasing and ${\bf e}$ a restricted growth sequence relative to ${\bf a}$, a weighted planar network all of whose weights are non-negative and whose path matrix is  $S^{{\bf a}, {\bf e}}$. We will achieve this construction in stages, first producing a weighted planar network whose path matrix is $S^{{\bf a}, {\bf e}}$ but which may have some negative weights, and then modifying the weight function in a way that makes all the negative weights non-negative, without changing the associated path matrix. 

\begin{conv}
Throughout the arguments that follow we will work with only one underlying planar network (the one shown in Figure \ref{fig1}), and from here on we will drop the qualifier ``planar''. All assignments of weights to this network will be of the form shown in (\ref{array-of-weights}), that is, we will always assign weight $1$ to the horizontal edges in the network. From here on, given an array of weights $V$ of the form shown in (\ref{array-of-weights}), rather than referring to ``the path matrix of the weighted network whose array of weights is $V$'', we will simply refer to ``the path matrix of the array of weights $V$''. 
\end{conv}

\begin{lemma} \label{lemma2-initial_network_construction}
For arbitrary ${\bf a}$ and ${\bf e}$, the following array of weights has $S^{{\bf a}, {\bf e}}$ as its path matrix:
\begin{equation} \label{orig-network}
\begin{array}{cccccccc}
a_1-e_1 & & & & & & \\
a_1-e_2 & a_2-e_1 & & & & & \\
a_1-e_3 & a_2-e_2 & a_3-e_1 & & & & \\
a_1-e_4 & a_2-e_3 & a_3-e_2 & a_4-e_1 & & & \\
a_1-e_5 & a_2-e_4 & a_3-e_3 & a_4-e_2 & a_5-e_1 & & \\
\vdots  & \vdots  & \vdots  & \vdots  &  & \ddots & \\
a_1-e_n & a_2-e_{n-1} & a_3-e_{n-2} & a_4-e_{n-3} & \cdots & a_{n-1}-e_2 & a_n-e_1 \\
\vdots  & \vdots & \vdots & \vdots &  & & & \ddots
\end{array}
\end{equation}
\end{lemma}

\medskip

\noindent {\em Proof}: Denote by $W$ the array shown in (\ref{orig-network}), and by $M$ its path matrix. Clearly the first row and column of $M$, as well as the main diagonal and everything above the main diagonal, agree with $S^{{\bf a}, {\bf e}}$, so we focus on the entries $M_{m,k}$ with $m > k \geq 1$.  

For each such fixed $m$ and $k$ the only weights which may appear on a path from $s_m$ to $t_k$ are those appearing in the subarray
\begin{equation} \label{square1}
\begin{array}{ccccc}
 a_1-e_{k+1} & a_2-e_k & a_3-e_{k-1} & \ldots & a_{k+1} - e_1 \\
 a_1-e_{k+2} & a_2-e_{k+1} & a_3-e_k & \ldots & a_{k+1} - e_2 \\
 a_1-e_{k+3} & a_2-e_{k+2} & a_3-e_{k+1} & \ldots & a_{k+1} - e_3 \\
 \vdots & \vdots & \vdots & & \vdots\\
 a_1-e_{m-1} & a_2 - e_{m-2} & a_3-e_{m-3} & \cdots & a_{k+1} - e_{m-k-1} \\
 a_1-e_m & a_2 - e_{m-1} & a_3-e_{m-2} & \cdots & a_{k+1} - e_{m-k}.
\end{array}
\end{equation}
(Specifically these are the weights in the $[i,j]$ position of $W$ for $k+1 \leq i \leq m$ and $1 \leq j \leq k+1$; there are $m-k$ rows and $k+1$ columns in the subarray). 

Now consider the array of weights $W'$ whose weights are identical to those of $W$, except that for $i=k+1, \ldots, m$ and $j=1, \ldots, k+1$ the weight $a_j-e_{i-j+1}$ is replaced by $a_j-e_{m-i+j}$. That is, $W'$ is obtained from $W$ by replacing the portion of $W$ shown in (\ref{square1}) with 
\begin{equation} \label{square2}
\begin{array}{ccccc}
 a_1-e_{m-k} & a_2-e_{m-k+1} & a_3-e_{m-k+2} & \ldots & a_{k+1} - e_m \\
 a_1-e_{m-k-1} & a_2-e_{m-k} & a_3-e_{m-k+1} & \ldots & a_{k+1} - e_{m-1} \\
 a_1-e_{m-k-2} & a_2-e_{m-k-1} & a_3-e_{m-k} & \ldots & a_{k+1} - e_{m-2} \\
 \vdots & \vdots & \vdots & & \vdots\\
 a_1-e_2 & a_2 - e_3 & a_3-e_4 & \cdots & a_{k+1} - e_{k+2} \\
 a_1-e_1 & a_2 - e_2 & a_3-e_3 & \cdots & a_{k+1} - e_{k+1}.
\end{array}
\end{equation}
Note that in going from (\ref{square1}) to (\ref{square2}) we are permuting ${\bf e}$ by the permutation that maps $e_i$ to $e_{m-i+1}$ for each $i=1, \ldots, m$.
 
A path from $s_m$ to $t_k$ in the underlying  network corresponds to a composition $b_1+\cdots+b_{k+1} = m-k$ of $m-k$ into $k+1$ non-negative parts, via: from $s_m$ take one horizontal step, then $b_1$ vertical steps along the first column of the network, then one horizontal step, then $b_2$ vertical steps along the second column of the network, and so on. 

In the  network whose array of weights is $W'$, the weight of the path corresponding to the composition $b_1+\cdots+b_{k+1} = m-k$ is a product of the form $\prod_{i=1}^{m-k} \left(a_{f(i)}-e_{g(i)}\right)$, where the sequence $(f(1), \ldots, f(m-k))$ consists of $b_1$ $1$'s, followed by $b_2$ $2$'s, and so on, and the sequence $(g(1), \ldots, g(m-k))$ starts $1, 2, \ldots, b_1$, then moves on to an increasing sequence of consecutive integers of length $b_2$ starting from $b_1+2$, and so on. In other words, the weight is
$$
\prod_{i=1}^{m-k} (a_{s_i-i+1} -e_{s_i})
$$
where $\{s_1, \ldots, s_{m-k}\}=\{1, \ldots, b_1, \widehat{b_1+1}, b_1+2,\ldots,b_1+b_2+1, \widehat{b_1+b_2+2},\ldots\}$ (the hats indicating missing elements). As $(b_1, \ldots, b_{k+1})$ runs over all $\binom{m}{m-k}$ compositions of $m-k$ into $k+1$ parts, the sets $\{s_1, \ldots, s_{m-k}\}$ run over all $\binom{m}{m-k}$ subsets of $\{1,\ldots,m\}$ of size $m-k$, and so we get that in the  network whose array of weights is $W'$ the sum of the weights of the paths from $s_m$ to $t_k$ is
\begin{equation} \label{int1}
\sum_{\displaystyle S=\{s_1,\ldots,s_{m-k}\} \subseteq \{1,\ldots,m\} \atop \displaystyle s_1<\ldots<s_{m-k}} \prod_{i=1}^{m-k} (a_{s_i-i+1} -e_{s_i}).
\end{equation}
By Lemma \ref{lem-gabes-formula} the expression in (\ref{int1}) is equal to $S^{{\bf a}, {\bf e}}(m,k)$. 
Thus by Lemma \ref{lem-sym} it is invariant under permutations of the $e_i$'s, $i=1, \ldots, m$, and so in particular, by considering the permutation that maps $e_i$ to $e_{m-i+1}$ for each $i=1, \ldots, m$ we get that
\begin{equation} \label{int2}
S^{{\bf a}, {\bf e}}(m,k) = \sum_{\displaystyle S=\{s_1,\ldots,s_{m-k}\} \subseteq \{1,\ldots,m\} \atop \displaystyle s_1<\ldots<s_{m-k}} \prod_{i=1}^{m-k} (a_{s_i-i+1} -e_{m-s_i+1}).
\end{equation}
Now using the same reasoning that led to (\ref{int1}) we see that in the  network whose array of weights is $W$ the sum of the weights of the paths from $s_m$ to $t_k$ is the right-hand side of (\ref{int2}), and the lemma is proved.
 \qed
 
\medskip

If $\inf {\bf a} \geq \sup {\bf e}$ then all the weights in $W$ are non-negative (monotonicity of ${\bf a}$ is not required for this), and so combining Lemma \ref{lemma2-initial_network_construction} and Lindstr\"om's Lemma we immediately get the following, discussed in the introduction.
\begin{cor} \label{cor-obsv-Gon}
For arbitrary ${\bf a}$ and ${\bf e}$ satisfying $\inf {\bf a} \geq \sup {\bf e}$, the matrix $S^{{\bf a}, {\bf e}}$ is totally non-negative.
\end{cor}

\medskip

Even if ${\bf a}$ is non-decreasing and ${\bf e}$ is a restricted growth sequence relative to ${\bf a}$, it may be that some of the weights in the array (\ref{orig-network}) are negative ($a_1-e_2$, for example).
We now describe a transformation that iteratively turns this array into one that has only non-negative weights, without changing the associated path matrix. 
 
\begin{lemma} \label{lemma3-transformation} 
For arbitrary ${\bf a}$ and ${\bf e}$, if $a_1=e_1$ then the following array of weights has $S^{{\bf a}, {\bf e}}$ as its path matrix:
\begin{equation} \label{mod-network}
\begin{array}{cccccccc}
a_1-e_1 & & & & & & \\
a_1-e_1 & a_2-e_2 & & & & & \\
a_1-e_1 & a_2-e_3 & a_3-e_2 & & & & \\
a_1-e_1 & a_2-e_4 & a_3-e_3 & a_4-e_2 & & & \\
a_1-e_1 & a_2-e_5 & a_3-e_4 & a_4-e_3 & a_5-e_2 & & \\
\vdots  & \vdots  & \vdots  & \vdots  &  & \ddots & \\
a_1-e_1 & a_2-e_n & a_3-e_{n-1} & a_4-e_{n-2} & \cdots & a_{n-1}-e_3 & a_n-e_2 \\
\vdots  & \vdots & \vdots & \vdots &  & & & \ddots
\end{array}
\end{equation}
\end{lemma}

\medskip

Note that the array shown in (\ref{mod-network}) is obtained from that shown in (\ref{orig-network}) by, in each row, moving the $-e_1$'s from the last position in the row to the first, and then shifting all other $-e_j$'s in the row one place to the right.

\medskip

\noindent {\em Proof (of Lemma \ref{lemma3-transformation})}: Denote by $W^{\rm p}$ the array shown in (\ref{mod-network}), and by $M^{\rm p}$ its path matrix. (We will shortly generalize the operation that transforms $W$ into $W^{\rm p}$, and refer to it as a ``pivoting'' operation; hence the notation $W^{\rm p}$ and $M^{\rm p}$.) Clearly the first row and column of $M^{\rm p}$, as well as the main diagonal and everything above the main diagonal, agree with $S^{{\bf a}, {\bf e}}$, so we focus on the entries $M^{\rm p}_{m,k}$ with $m > k \geq 1$.  

For each such fixed $m$ and $k$ the only weights which may appear on a path from $s_m$ to $t_k$ are a subset of those those appearing in the subarray consisting of the weights at the $[i,j]$ position for $m \geq i \geq j \geq 1$. Because $a_1=e_1$, this subarray takes the following form:
$$
\begin{array}{ccccccc}
0 & & & & & & \\
0 &a_2-e_2 & & & & & \\
0 & a_2-e_3 & a_3-e_2 & & & & \\
0 &a_2-e_4 & a_3-e_3 & a_4-e_2 & & & \\
0 & a_2-e_5 & a_3-e_4 & a_4-e_3 & a_5-e_2 & & \\
\vdots  & \vdots  & \vdots  & \vdots  &  & \ddots & \\
0 & a_2-e_m & a_3-e_{m-1} & a_4-e_{m-2} & \cdots & a_{m-1}-e_3 & a_m-e_2.
\end{array}
$$
The sum of the weights of the paths from $s_m$ to $t_k$ in the  network whose array of weights is $W^{\rm p}$ is evidently the same as the sum of the weights of the paths from $s_{m-1}$ to $t_{k-1}$ in the  network whose array of weights has the following as its first $m-1$ rows:
$$
\begin{array}{cccccc}
a_2-e_2 & & & & & \\
a_2-e_3 & a_3-e_2 & & & & \\
a_2-e_4 & a_3-e_3 & a_4-e_2 & & & \\
a_2-e_5 & a_3-e_4 & a_4-e_3 & a_5-e_2 & & \\
\vdots  & \vdots  & \vdots  &  & \ddots & \\
a_2-e_m & a_3-e_{m-1} & a_4-e_{m-2} & \cdots & a_{m-1}-e_3 & a_m-e_2.
\end{array}
$$ 
From the proof of Lemma \ref{lemma2-initial_network_construction} this quantity is symmetric in $e_2, \ldots, e_m$. So, if $W^{\rm p'}$ is the array of weights obtained from $W^{\rm p}$ by the transformation $e_2\rightarrow e_m$, $e_3 \rightarrow e_{m-1}$, et cetera (that is, by replacing $e_i$ with $e_{m-i+2}$ for $i=2, \ldots, m$), then although this perhaps changes the path matrix, it does not change the sum of the weights of the paths from $s_m$ to $t_k$ (that is, the $(m,k)$ entry of the path matrix). 

Now consider the array of weights $W^{\rm p''}$ obtained from $W^{\rm p'}$ by changing the weight at the $[i,1]$ position from $a_1-e_1$ to $a_1 - e_{m-i+1}$, for $i=1, \ldots, m$. The first $m$ rows of $W^{\rm p''}$ have the following form:
$$
\begin{array}{ccccccc}
a_1-e_m & & & & & & \\
a_1 -e_{m-1} &a_2-e_m & & & & & \\
a_1 -e_{m-2} & a_2-e_{m-1} & a_3-e_m & & & & \\
a_1 -e_{m-3} &a_2-e_{m-2} & a_3-e_{m-1} & a_4-e_m & & & \\
a_1 -e_{m-4} & a_2-e_{m-3} & a_3-e_{m-2} & a_4-e_{m-1} & a_5-e_m & & \\
\vdots  & \vdots  & \vdots  & \vdots  &  & \ddots & \\
a_1 -e_1 & a_2-e_2 & a_3-e_3 & a_4-e_4 & \cdots & a_{m-1}-e_{m-1} & a_m-e_m.
\end{array}
$$
The sum of the weights of the paths from $s_m$ to $t_k$ in the network whose array of weights is $W^{\rm p}$ is, as has been observed, the same as that for the network whose array of weights is $W^{\rm p'}$. We now argue that this sum is the same as that for the network whose array of weights is $W^{\rm p''}$. Indeed, the only weights that have (potentially) changed in going from $W^{\rm p'}$ to $W^{\rm p''}$ are those in the $[1,1]$ through $[m-1,1]$ positions, and any path from $s_m$ to $t_k$ that uses the $[k,1]$ edge for some $k < m$ must also use the $[m,1]$ edge, which has weight $0$.

The only weights in $W^{\rm p''}$ which may appear on a path from $s_m$ to $t_k$ are those in the $[i,j]$ position of $W^{\rm p''}$ for $k+1 \leq i \leq m$ and $1 \leq j \leq k+1$. The rectangular subarray of weights in $W^{\rm p''}$ in those positions is exactly the subarray (\ref{square2}), and so the proof of Lemma \ref{lemma2-initial_network_construction} shows that in the network whose array of weights is $W^{\rm p''}$, the sum of the weights of the paths from $s_m$ to $t_k$ is $S^{\bf{a}, {\bf e}}(m,k)$.
\qed

\medskip

We refer to the operation that transforms the array of weights $W$ of Lemma \ref{lemma2-initial_network_construction} (shown in (\ref{orig-network})) to the array of weights $W^{\rm p}$ of Lemma \ref{lemma3-transformation} (shown in (\ref{mod-network})) as {\em pivoting} on the $[1,1]$ position of the array. We now define a more general pivoting operation:
\begin{defn}
Let $V$ be a doubly infinite lower triangular array of weights (as shown in (\ref{array-of-weights})) associated with the  network shown in Figure \ref{fig1}. If for each $m \geq 1$ and $1 \leq k \leq m$ the weight in the $[m,k]$ position of $V$ is of the form $a_{f(m,k)}-e_{g(m,k)}$ (for some functions $f, g$) then we denote by $V^{[m,k]}$ the array of weights constructed from $V$ by the following process:
\begin{itemize}
\item the weight in the $[m,k]$ position remains unchanged; 
\item in row $m+1$, the weights $a_{f(m+1,k)}-e_{g(m+1,k)}$ and $a_{f(m+1,k+1)}-e_{g(m+1,k+1)}$ (in the $[m+1,k]$ and $[m+1,k+1]$ positions, respectively) are replaced with $a_{f(m+1,k)}-e_{g(m+1,k+1)}$ and $a_{f(m+1,k+1)}-e_{g(m+1,k)}$;
\item in general, for $\ell \geq 1$ the weights
$$
a_{f(m+\ell,k)}-e_{g(m+\ell,k)},\ a_{f(m+\ell,k+1)}-e_{g(m+\ell,k+1)},\ \ldots,\ a_{f(m+\ell,k+\ell)}-e_{g(m+\ell,k+\ell)} 
$$  
(in the $[m+\ell,k]$ through $[m+\ell,k+\ell]$ positions, respectively) are replaced with
$$
a_{f(m+\ell,k)}-e_{g(m+\ell,k+\ell)},\ a_{f(m+\ell,k+\ell)}-e_{g(m+\ell,k)},\ \ldots,\ a_{f(m+\ell,k+\ell)}-e_{g(m+\ell,k+\ell-1)};  
$$
\item and all other weights remain unchanged.
\end{itemize}
We refer to $V^{[m,k]}$ as the array of weights obtained from $V$ by
{\em pivoting} on the $[m,k]$ position.
\end{defn}

For a doubly infinite lower triangular array (such as the one shown in (\ref{array-of-weights})) we refer to the triangle consisting of the $[m+\ell_1,k+\ell_2]$ positions for all $\ell_1 \geq 0$ and $0 \leq \ell_2 \leq \ell_1$ as the triangle {\em headed} at the $[m,k]$ position, and we refer to the collection of positions that are in the $i$th column of the array, for $i \geq k$, but that are not in the triangle headed at the $[m,k]$ position, as the positions that lie {\em above} the triangle. Figure \ref{fig2} shows a portion of the triangle headed at the $[3,2]$ position (the bolded entries), and the positions lying above that triangle (the italicized entries).

\begin{figure}[ht!]
\begin{center}
$$
\begin{array}{ccccccccc}
 & [1,1] & & & & & & & \\
& [2,1] & \emph{[2,2]} & & & & & & \\
 & [3,1] & {\bf [3,2]} & \emph{[3,3]} & & & & & \\
 & [4,1] & {\bf [4,2]} & {\bf [4,3]} & \emph{[4,4]} & & & & \\
 & [5,1] & {\bf [5,2]} & {\bf [5,3]} & {\bf [5,4]} & \emph{[5,5]} & & & \\
 & \vdots & \vdots &  &  & \ddots & \ddots  & & \\
 & [m,1] & {\bf [m,2]} & {\bf [m,3]} & {\bf [m,4]} & {\bf \cdots} & {\bf [m,m-1]} & \emph{[m,m]} &  \\ 
 & \vdots & \vdots &  &  &  & ~~~~~ \ddots &  \ddots & 
\end{array}
$$
\caption{The triangle headed at the $[3,2]$ position (bolded entries), and the positions that lie above the triangle (italicized entries).} \label{fig2}
\end{center}
\end{figure}

We now generalize Lemma \ref{lemma3-transformation}.
\begin{lemma} \label{lema3-transformation2}
Let ${\bf a}$ and ${\bf e}$ be arbitrary. Let $W$ be the array of weights shown in (\ref{orig-network}). If the weight in the $[m,k]$ position of $W$ is $0$, then the path matrix of the array of weights $W^{[m,k]}$ is the same as that of $W$, that is, it is $S^{{\bf a}, {\bf e}}$.

Furthermore, let $(m_1, m_2, \ldots)$ and $(k_1, k_2, \ldots)$ be sequences satisfying that for each $i\geq 1$, the $[m_{i+1},k_{i+1}]$ position is located in the triangle headed at the $[m_i,k_i]$ position. Let $\overline{W}$ be the array of weights obtained from $W$ by first pivoting on the $[m_1,k_1]$ position, then pivoting on the $[m_2,k_2]$ position of the resulting array, and so on. If the weight at each position at which pivoting occurs is $0$ (at the moment when the pivoting occurs at that position), then the path matrix of $\overline{W}$ is still $S^{{\bf a}, {\bf e}}$.    
\end{lemma}

\medskip

\noindent {\em Proof}:
We begin with the first statement. Say that an edge in the  network is {\em in} the triangle headed at the $[m,k]$ position if it is the $[i,j]$ edge of the network for some $i, j$ such that the $[i,j]$ position is in the triangle (see the paragraph after (\ref{array-of-weights}) for relevant definitions). Fix a source $s_p$ and sink $t_q$. Each path from $s_p$ to $t_q$ starts with a (non-empty) path $A$ consisting of edges all not in the triangle headed at the $[m,k]$ position, then continues with a (possibly empty) path $B$ consisting of edges all in the triangle, and then ends with a (also possibly empty) path $C$ consisting of edges all not in the triangle. The collection of paths from $s_p$ to $t_q$ thus can be partitioned into a collection of blocks indexed by pairs $(v_f,v_\ell)$, where $v_f$ is the first vertex along $B$ and $v_\ell$ is the last vertex along $B$, and an exceptional block consisting of those paths for which $B$ (and so also $C$) is empty. 

For a block indexed by the pair $(v_f,v_\ell)$, the sum of the weights of the paths from $s_p$ to $t_q$ is the product of three factors: the sum of the weights of the paths from $s_p$ to $v_f$, the sum of the weights of the paths from $v_f$ to $v_l$, and sum of the weights of the paths from $v_l$ to $t_q$. The first and third of these sums remain unchanged after pivoting on the $[m,k]$ position, because the pivoting does not change the weight at any position not in the triangle headed at the $[m,k]$ position. The middle sum also remains unchanged after pivoting, by Lemma \ref{lemma3-transformation} (applied in the obvious way to the array of weights in the triangle headed at the $[m,k]$ position). For the exceptional block, the sum of the weights of the paths from $s_p$ to $t_q$ clearly also remains unchanged after pivoting on the $[m,k]$ position. Summing over blocks, the first statement of the lemma follows.     

The second statement of the lemma is obtained by iterating the above argument.
\qed

\medskip 

We can now fairly swiftly present the proof of Theorem \ref{thm-main}.

\medskip

\noindent {\em Proof (of Theorem \ref{thm-main})}: Let ${\bf a}$ be non-decreasing. We begin by arguing that if ${\bf e}$ is a restricted growth sequence relative to ${\bf a}$, then $S^{{\bf a}, {\bf e}}$ is totally non-negative (item \ref{mainthm-item1}).

\begin{itemize}

\item If all $e_i$ are at most $a_1$, then the array of weights $W$ shown in (\ref{orig-network}) evidently has all non-negative weights, and by Lemma \ref{lemma2-initial_network_construction} has path matrix $S^{{\bf a}, {\bf e}}$. By Lemma \ref{lem-tp} (Lindstr\"om's Lemma) we are done.

\item If it is not the case that all $e_i$ are at most $a_1$, then there is some index $j$ such that $e_j=a_1$ and $e_{j'}<a_1$ for all $j'<j$. We pivot on the $[j,1]$ position of $W$. Note that the weight in this position is $a_1-e_j=0$, so from the first part of Lemma \ref{lema3-transformation2} the path matrix of the resulting array of weights $W^{[j,1]}$ is $S^{{\bf a}, {\bf e}}$. Notice that all weights in the first column of $W^{[j,1]}$ are either positive (the weights in the first $j-1$ rows) or $0$ (the remaining weights), and that all weights in $W^{[j,1]}$ that lie above the triangle headed at the $[j,1]$ position are positive (they are positive in $W$ --- here we use that ${\bf a}$ is non-decreasing --- and remain unchanged after pivoting). In other words, after pivoting all weights in the new array in positions outside the triangle headed at the $[j,1]$ position are non-negative.

\item If all $e_i$ for $i > j$ are at most $a_2$, then array $W^{[j,1]}$ has only non-negative weights, and again by Lemma \ref{lem-tp} we are done. If not, there is some index $j'$ such that $e_{j'}=a_2$ and $e_{j''}<a_2$ for all $j <j''<j'$. We now pivot on the $[j',2]$ position in $W^{[j,1]}$ (which has weight $a_2-e_{j'}=0$). Because the $[j',2]$ position is in the triangle headed at the $[j,1]$ position, we can apply the second part of Lemma \ref{lema3-transformation2} to conclude that the path matrix of the resulting array of weights is still $S^{{\bf a}, {\bf e}}$. Arguing as before, the new array of weights has non-negative weights outside the triangle headed at the $[j',2]$ position.  

\item Iterating this process (either finitely many times or countably many times, depending on whether ${\bf a}$ and ${\bf e}$ are finite or countably infinite) we arrive at an array of weights all of whose entries are non-negative and whose path matrix is $S^{{\bf a}, {\bf e}}$; the result now follows from Lemma \ref{lem-tp}. 

\end{itemize}

To complete the proof of Theorem \ref{thm-main}, we show that if ${\bf e}$ is not a restricted growth sequence relative to ${\bf a}$, then $S^{{\bf a}, {\bf e}}$ is {\em not} totally non-negative, and that moreover the failure of total non-negativity is witnessed by a negative matrix entry (item \ref{mainthm-item2}).

\begin{itemize}

\item Suppose that the failure of ${\bf e}$ to be a restricted growth sequence relative to ${\bf a}$ is witnessed by some index $j$ such that $e_i < a_1$ for all $i < j$, and $e_j > a_1$. Then evidently the path matrix of the array of weights $W$ has the negative entry $(a_1-e_j)(a_1-e_{j-1})\cdots (a_1-e_1)$ --- it is the $(j,0)$ entry. 

\item Otherwise, there is some index $j$ such that $e_i < a_1$ for all $i < j$, and $e_j = a_1$. Consider the array of weights $W^{[j,1]}$. As established in the proof of item \ref{mainthm-item1} above, $W^{[j,1]}$ has path matrix $S^{{\bf a}, {\bf e}}$. Also, it has strictly positive weights in the first $j-1$ entries of the first column, the weights in the rest of the first column are all $0$, and all weights above the triangle headed at the $[j,1]$ position are strictly positive. 

Now suppose that the failure of ${\bf e}$ to be a restricted growth sequence relative to ${\bf a}$ is witnessed by some index $j'$ such that $e_i < a_2$ for all $j <i < j'$, and $e_{j'} > a_2$. Evidently the $(j',1)$ entry of the path matrix of the array $W^{[j,1]}$ is negative, because all paths from $s_{j'}$ to $t_1$ that do not have weight $0$ have a weight which is a product of strictly positive terms, together with the term $a_2-e_{j'}$, which is negative.

\item Repeating this argument, we obtain the general result that if the earliest witness of the failure of ${\bf e}$ to be a restricted growth sequence relative to ${\bf a}$ is some index $\tilde{j}$ with $e_{\tilde{j}}>a_\ell$ for some $\ell$, then the $(\tilde{j},\ell-1)$ entry of $S^{{\bf a}, {\bf e}}$ is negative.

\end{itemize}

\qed

\subsection*{Acknowledgements}

We are grateful to Gabriel Conant for conjecturing (\ref{formula-Gabe2}), to K. Gonzales for pointing out some useful references, and to two referees for their detailed and helpful comments.

\end{document}